\documentclass[preprint,11pt]{elsarticle}

\newtheorem{theorem}{Theorem}[section]
\newtheorem{lemma}{Lemma}[section]
\newtheorem{corollary}{Corollary}[section]

\newtheorem{remark}{Remark}[section]
\newcommand{\ignore}[1]{}{}

\usepackage{amssymb}
\usepackage{amsmath}

\usepackage{color}
\usepackage{soul}
\usepackage[dvipsnames]{xcolor}

\usepackage[colorlinks=true, linkcolor=blue, citecolor=blue]{hyperref}

\def\1{{{\mbox{${\rm{1\negthinspace\negthinspace I}}$}}}}

\newcommand\beq{\begin{equation}}
\newcommand\eeq{\end{equation}}

\usepackage{ifthen}
\usepackage{xkeyval}
\usepackage{todonotes}
\begin{document}
\begin{frontmatter}
\title{Cram\'{e}r-type moderate deviations for Euler-Maruyama scheme for SDE}
\author[cor1]{Xiequan Fan}
\author[cor2]{Haijuan Hu}
\author[cor3]{Lihu Xu}
\address[cor1]{Center for Applied Mathematics,
Tianjin University, Tianjin 300072,  China;}
\address[cor2]{School of Mathematics and Statistics, Northeastern University at Qinhuangdao, Qinhuangdao, China;}
\address[cor3]{Department of Mathematics, Faculty of Science and Technology, University of Macau, Macau, China.}

\begin{abstract}
In this paper, we establish normalized and self-normalized Cram\'{e}r-type moderate deviations for Euler-Maruyama scheme for SDE.
As a consequence of our results,   Berry-Esseen's bounds and moderate deviation principles are also obtained.
Our normalized Cram\'{e}r-type moderate deviations refines the recent work of [Lu, J., Tan, Y., Xu, L., 2022. Central limit theorem and self-normalized Cram\'{e}r-type moderate deviation for Euler-Maruyama scheme. \emph{Bernoulli} \textbf{28}(2): 937--964].
\end{abstract}
\begin{keyword}  Euler-Maruyama scheme;  Cram\'{e}r-type moderate deviations;  Self-normalized sequences; Berry-Esseen's bounds
\vspace{0.3cm}
\MSC primary  60F10; 60E05;  secondary  60H10;  62F12
\end{keyword}
\end{frontmatter}




\section{Introduction}
Consider the following stochastic differential equation (SDE) on $\mathbf{R}^d$:
\begin{eqnarray}\label{sde}
 dX_t=g(X_t)dt + \sigma dB_t,\ \ \ X_0=\mathbf{x}_0,
\end{eqnarray}
where $B_t$ is a $d$-dimential standard Brownian motion,   $\sigma$ is  an invertible $d\times d$ matrix  and   $g: \mathbf{R}^d \rightarrow \mathbf{R}^{d}$ satisfies the following assumption. There exist  constants $L, K_1 >0$ and $K_2\geq 0$ such that for every $\mathbf{x}, \mathbf{y} \in \mathbf{R}^d$,
 \begin{eqnarray*}
\|g(\mathbf{x}) -g (\mathbf{y})\|_2 \!\!\!& \leq&\!\!\! L \|\mathbf{x}-\mathbf{y}\|_2, \\
\langle g(\mathbf{x})-g(\mathbf{y}), \mathbf{x}-\mathbf{y}   \rangle \!\!\!&\leq&\!\!\! - K_1 \|\mathbf{x}-\mathbf{y}\|_2 + K_2,
\end{eqnarray*}
where $\langle \mathbf{x},   \mathbf{y}  \rangle$ stands for inner product of $\mathbf{x}$ and $\mathbf{y}$, and $\| \cdot\|_2$ is the Euclidean norm.
Moreover, assume that $g(\mathbf{x})$ is second order differentiable and the  second order derivative of $g$ is bounded.
Given step size $\eta \in (0, 1),$ the Euler-Maruyama scheme for SDE (\ref{sde}) is given by
 \begin{eqnarray}\label{fsdf3}
  \theta_{k+1} = \theta_k + \eta g(\theta_k) + \sqrt{\eta} \sigma \xi_{k+1},\ \ \ \ k \geq 0,
  \end{eqnarray}
where   $(\xi_{k})_{k\geq1}$ are independent and identically distributed (i.i.d.)\ standard $d$-dimensional normal random vectors.
It is known that SDE (\ref{sde}) and $(\theta_k)_{k\geq0}$ are both ergodic and admit invariant measures, denoted by  $\pi$ and $\pi_\eta$ respectively; see Lemma 2.3
of Lu, Tan and Xu \cite{LTX22}. Moreover,  when $\sigma$ is the  identity matrix and $g(\mathbf{x})$  is third order differentiable with
an appropriate growth condition,  Fang, Shao and Xu \cite{FSX19} have proved  that the Wasserstein-1 distance between $\pi$ and $\pi_\eta$
is in order of  $\,\eta^{1/2}$, up to a logarithmic correction.

 Denote $ C_{b}^2(\mathbf{R}^d, \mathbf{R})$
the collection of all bounded $2$-th order continuously differentiable functions. Given an $h \in C_{b}^2(\mathbf{R}^d, \mathbf{R})$,  denote $\varphi$ the solution to the following Stein's equation:
\begin{eqnarray}\label{stein's}
 h-\pi(h)=\mathcal{A} \varphi,
\end{eqnarray}
where $\mathcal{A} $ is the generator of SDE (\ref{sde}) defined as follows:
$$  \mathcal{A} \varphi(\mathbf{x})= \langle g(\mathbf{x}), \nabla  \varphi(\mathbf{x})\rangle + \frac12 \langle \sigma\sigma^T, \nabla^2 \varphi(\mathbf{x})   \rangle_{\textrm{HS}} ,$$
with $T$  the transport operator and $\langle A, B  \rangle_{\textrm{HS}}:=\sum_{i=1}^d\sum_{j=1}^d a_{ij} b_{ij} $ for
$A=(a_{ij})_{d\times d},  B=(b_{ij})_{d\times d} \in \mathbf{R}^{d\times d}.$
For a small $\eta \in (0, 1),$ define
$$  \Pi_\eta(\cdot)=  \frac{1}{[\eta^{-2}]} \sum_{k=0}^{[\eta^{-2}]-1  } \delta_{\theta_k}(\cdot) ,$$
where $\delta_{y}(\cdot)$ is the Dirac measure of $y.$ Then $ \Pi_\eta$ is an asymptotically consistent statistic of $\pi $  as $\eta \rightarrow 0.$
 We also denote
$$W_\eta=\frac{    \eta^{-1/2}   (  \Pi_\eta(h) - \pi (h) )}{ \sqrt{ \mathcal{Y}_{\eta}}} \ \ \ \ \ \textrm{with} \ \ \ \  \mathcal{Y}_\eta=\frac{1}{[\eta^{-2}]} \sum_{k=0}^{[\eta^{-2}]-1  }\|\sigma^T \nabla \varphi(\theta_k)\|_2^2 . $$
 Recently,
Lu, Tan and Xu \cite{LTX22}  proved the following  normalized Cram\'{e}r-type moderate deviation. If $\theta_0 \sim \pi_\eta$ and $h \in C_{b}^2(\mathbf{R}^d, \mathbf{R})$, then
\begin{eqnarray}\label{ffgds}
  \frac{\mathbf{P} \big(  W_\eta > x   \big)}{1- \Phi(x)}   =1+O\big((1+x)  \eta^{1/6} \big)
\end{eqnarray}
holds uniformly for $ c \eta^{1/6} \leq x =o(\eta^{-1/6}) $ as $\eta \rightarrow 0.$
In this paper, we  give an improvement on (\ref{ffgds}). In particular, our result implies that
\begin{eqnarray}\label{dsdgs}
 \frac{\mathbf{P} \big(  W_\eta > x   \big)}{1- \Phi(x)} =1+o\bigg(   x^2 \eta^{1/2}   + (1+ x) (\eta  |\ln \eta|)^{1/2}      \bigg)
\end{eqnarray}
holds uniformly for $ 0 \leq x =o(\eta^{-1/4}) $ as $\eta \rightarrow 0.$
Compared to  (\ref{ffgds}),  equality (\ref{dsdgs}) holds for a much larger range. Moreover, from (\ref{dsdgs}), we obtain the following Berry-Esseen bound
\begin{eqnarray}\label{dsdsfgs}
 \sup_{x \in \mathbf{R}}\Big|\mathbf{P}( W_\eta \leq  x  ) - \Phi \left(  x\right)\Big|  \leq c  \, (\eta  |\ln \eta|)^{1/2}.
\end{eqnarray}
Notice that   the limit $\lim_{\eta \rightarrow 0}  \eta^{-1/2}   (  \Pi_\eta(h) - \pi (h) ) $ has a normal distribution. Thus the best possible convergence rate of Berry-Esseen's  bound  is in order of $\eta^{1/2}$.
 Thus the convergence rate  in the last Berry-Esseen  bound (\ref{dsdsfgs}) is close to the best possible one $ \eta^{1/2}$,   up to a logarithmic correction $|\ln \eta|^{1/2}$.
In particular, we further  establish the following  self-normalized Cram\'{e}r-type moderate deviation.
Denote
\begin{eqnarray}\label{ddfsn}
   S_\eta=\frac{    \eta^{-1/2}   (  \Pi_\eta(h) - \pi (h) )}{ \sqrt{\mathcal{V}_\eta}} \ \ \   \textrm{with}  \ \ \   \mathcal{V}_\eta=\frac{1}{ \eta\, [\eta^{-2}] } \sum_{k=0}^{[\eta^{-2}]-1  } \big\|( \theta_{k+1} - \theta_k - \eta g(\theta_k) )^T \nabla \varphi (\theta_k) \big\|_2^2.
 \end{eqnarray}
We also show that (\ref{dsdgs}) and (\ref{dsdsfgs}) hold also when $W_\eta$ is replaced by $S_\eta$.
As $(\theta_k)_{0\leq k  \leq [\eta^{-2}]-1}$  are observable, then  $S_\eta$ is a self-normalized process.
The moderate deviation expansion with respect to $S_\eta$  is called as
self-normalized Cram\'{e}r-type moderate deviation. Self-normalized Cram\'{e}r-type moderate deviation plays an important role
 in statistical inference of $\pi (h)$, because in practice one usually does not know the exact values of the matrix  $\sigma$ and the factor $\mathcal{V}_\eta$ does not depend  on the invertible matrix  $\sigma.$

Throughout the paper, denote $C$ a positive constant,  and denote $c$ a positive constant depending only on $ L, K_1, K_2, g, \|g(0)\|_2 $ and $\sigma$. The  exact values of $C$ and $c$ may vary from line to line.
All over the paper, $|| \cdot  ||$ stands for the Euclidean norm for higher rank tensors.

\section{Main results}
\setcounter{equation}{0}

We have the following  normalized Cram\'{e}r-type moderate deviation.
\begin{theorem}\label{the3.1}
Let $\theta_0 \sim \pi_\eta$ and $h \in C_{b}^2(\mathbf{R}^d, \mathbf{R})$. Then the following  inequality
\begin{eqnarray*}
\bigg| \ln \frac{\mathbf{P} \big(  W_\eta > x   \big)}{1- \Phi(x)} \bigg|  \leq c  \bigg(      x^3  \eta  +x^2 \eta^{1/2}   + (1+ x) (\eta  |\ln \eta|)^{1/2}      \bigg)
\end{eqnarray*}
holds uniformly for $0 \leq x \leq \eta^{-3/4} $. In particular, it implies that
\begin{eqnarray*}
 \frac{\mathbf{P} \big(  W_\eta > x   \big)}{1- \Phi(x)} =1+o\bigg(   x^2 \eta^{1/2}   + (1+ x) (\eta  |\ln \eta|)^{1/2}      \bigg)
\end{eqnarray*}
holds uniformly for $0 \leq x =o  (  \eta^{-1/4}   ) $ as $\eta \rightarrow 0.$ Moreover, the same results hold when $W_\eta $ is replaced by $-W_\eta$.
\end{theorem}

 From Theorem \ref{the3.1},    we have the following Berry-Esseen bound for $W_\eta$.
\begin{corollary}\label{coro1}
Assume the conditions of Theorem \ref{the3.1}. Then  the following inequality holds
\begin{equation}
 \sup_{x \in \mathbf{R}}\Big|\mathbf{P}( W_\eta \leq  x  ) - \Phi \left(  x\right)\Big|  \leq c  \, (\eta  |\ln \eta|)^{1/2}  .
\end{equation}
\end{corollary}

Notice that from (\ref{ffgds}), one may obtain a Berry-Esseen's bound of order $ \eta  ^{1/6},$ which is slower the one in Corollary \ref{coro1}.
Moreover, the convergence rate of Berry-Esseen's bound in Corollary \ref{coro1} is close to the convergence rate in the
 Wasserstein-1 distance between $\pi$ and $\pi_\eta$, which is of order $\eta^{1/2}$ up to a logarithmic correction.

From Theorem \ref{the3.1},
by an argument similar to the proof of Corollary 2.2 in \cite{FGLS19}, we easily obtain the following moderate deviation principle  (MDP) result.
\begin{corollary}\label{coro2}
Assume the conditions of Theorem \ref{the3.1}.
Let $(a_\eta)$ be any sequence of real numbers satisfying $a_\eta \rightarrow \infty$ and $\eta^{3/4} a_\eta   \rightarrow 0$
as $\eta \rightarrow 0$.  Then  for each Borel set $B$,
\begin{eqnarray}
- \inf_{x \in B^o}\frac{x^2}{2} \!\!\!&\leq &\!\!\! \liminf_{\eta \rightarrow 0}\frac{1}{a_\eta^2}\ln \mathbf{P}\bigg(\frac{  W_\eta}{a_\eta }  \in B \bigg) \nonumber \\
 &\leq&\!\!\! \limsup_{\eta\rightarrow 0}\frac{1}{a_\eta^2}\ln \mathbf{P}\bigg(\frac{ W_\eta}{a_\eta }    \in B \bigg) \leq  - \inf_{x \in \overline{B}}\frac{x^2}{2} \, ,   \label{MDP}
\end{eqnarray}
where $B^o$ and $\overline{B}$ denote the interior and the closure of $B$, respectively.
\end{corollary}

Self-normalized limit theory for independent random variables has been studied in depth  in the past twenty-five years. See, for instance, Shao \cite{S97}
for self-normalized large deviations,  and Jing, Shao and Wang \cite{JSW03} and Shao and Zhou \cite{SZ16} for self-normalized Cram\'{e}r-type moderate  deviations.
A few results for dependent random variables, we refer to Chen, Shao,   Wu and Xu \cite{CSWX16}  and  \cite{FGLS19}.
In the next theorem, we present a self-normalized Cram\'{e}r-type moderate deviation for the Euler-Maruyama scheme for SDE (\ref{sde}).
Recall the definition of $S_\eta$ in (\ref{ddfsn}).
We have the following moderate deviation expansions.
\begin{theorem}\label{the3.3}
Let $\theta_0 \sim \pi_\eta$ and $h \in C_{b}^2(\mathbf{R}^d, \mathbf{R})$. Then the following  inequality
\begin{eqnarray*}
\bigg| \ln \frac{\mathbf{P} \big(  S_\eta > x   \big)}{1- \Phi(x)} \bigg|  \leq c  \bigg(      x^3  \eta  +x^2 \eta^{1/2}   + (1+ x) (\eta  |\ln \eta|)^{1/2}      \bigg)
\end{eqnarray*}
holds uniformly for $0 \leq x \leq \eta^{-3/4} $. In particular, it implies that
\begin{eqnarray*}
 \frac{\mathbf{P} \big(  S_\eta > x   \big)}{1- \Phi(x)} =1+o\bigg(   x^2 \eta^{1/2}   + (1+ x) (\eta  |\ln \eta|)^{1/2}      \bigg)
\end{eqnarray*}
holds uniformly for $0 \leq x =o  (  \eta^{-1/4}   ) $ as $\eta \rightarrow 0.$ Moreover, the same results hold when $S_\eta $ is replaced by $-S_\eta$.
\end{theorem}

We  call the results in Theorem \ref{the3.3} as self-normalized Cram\'{e}r-type moderate deviations,
because the normalized factor $\mathcal{V}_\eta$ does  not depend  on the invertible matrix  $\sigma.$
Such type results play  an important role in statistical inference of $\pi (h)$, since
in practice one usually does not know the exact values of the matrix  $\sigma.$

By arguments similar to the proofs of Corollaries \ref{coro1} and \ref{coro2}, from Theorem \ref{the3.3}, it is easy to see that the assertions in
Corollaries \ref{coro1} and \ref{coro2}  remain valid when $W_\eta$ is replaced by $S_\eta$.

\begin{remark}
The assumption $\theta_0 \sim \pi_\eta$ in Theorems \ref{the3.1} and  \ref{the3.3} is not essential. Thanks to the exponential ergodicity of the EM scheme, one can extend Theorems \ref{the3.1} and  \ref{the3.3} to the case in which $\theta_0$ is subgaussian distributed. 
Indeed, from the proof of (3.2) in Lu, Tan and Xu \cite{LTX22},
one can see that Lemma  3.2   holds also when $\theta_0$ is subgaussian distributed. The advantage of taking $\theta_0 \sim \pi_\eta$ is that in their calculations, the terms describing the difference between the distribution of $\theta_{k}$ and $\pi_\eta$ will vanish, while in general case, one has to use exponential ergodicity of $\theta_k$ to bound the difference. Since $\theta_k$ converges to $\pi_\eta$ exponentially fast, the difference will not put an essential difficulty.
In our paper, we assumed the same condition as in Lu, Tan and Xu \cite{LTX22} and thus only considered the case of $\theta_0 \sim \pi_\eta$.
\end{remark}

\section{Preliminary lemmas}\label{sec5}
\setcounter{equation}{0}
In the proof of Theorem \ref{the3.1}, we   need the following three lemmas of  Lu, Tan and Xu \cite{LTX22}, see Lemmas 3.1, 3.3 and 5.1 therein.
\begin{lemma}\label{lm22}
 Let   $h \in C_{b}^2(\mathbf{R}^d, \mathbf{R})$. Then
 $$ || \bigtriangledown^k \varphi  || \leq c,\ \ \ k=0,1,2,3,4,$$
  where   $c $ depends  on $g$ and $\sigma.$
\end{lemma}

\begin{lemma}\label{lm23}
If $\theta_0 \sim \pi_\eta$, then there exists a constant $\gamma >0$, depending on $L, K_1, K_2, \|g(0)\|_2 $ and $\sigma$, such that
 $$\mathbf{E}\bigg( \exp\Big\{\gamma \eta  \sum_{k=0}^{m-1}   \|     g(\theta_k)  \|_2^2 \Big \} \bigg) \leq   c_1 e^{c_2 \eta^{-1}}  $$
 and for all $x>0,$
  $$\mathbf{P}\bigg( \eta  \sum_{k=0}^{m-1}   \|     g(\theta_k)  \|_2^2 >x \bigg) \leq c_1 e^{c_2 \eta^{-1}}e^{-c_3 x} .$$
\end{lemma}

\begin{lemma}\label{lm24}
For any $k \in \mathbf{N} $,  it holds for all $y>0,$
  $$\mathbf{P}\bigg(    \Big |  \sum_{i=0}^{k-1  }\|\sigma^T \nabla \varphi(\theta_i)\|_2^2 - k \pi_{\eta} (\|\sigma^T \nabla \varphi \|_2^2)    \Big|  > y \bigg) \leq 2\, e^{- c \, y^2 k^{-1}} ,\ \ \ \,$$
  where $c$ depends on $g$ and $\sigma.$
\end{lemma}

In the proof of Theorem \ref{the3.1}, we also make use of the following lemma.

\begin{lemma}\label{lm21}
Let $(\zeta _i,\mathcal{F}_i)_{i\geq1}$ be a sequence of martingale differences.
Assume there exist  positive constants $c$  and $\alpha \in (0,1]$ such that
\begin{eqnarray}\label{fcond14f}
1\leq u_n:=\sum_{i=1}^n \Big\|\mathbf{E} (\zeta_i^2 \exp\{ c |\zeta_{i}|^{\alpha} \} |\mathcal{F}_{i-1}) \Big\|_\infty  < \infty.
\end{eqnarray}
Then there exits  a positive constant $c_\alpha$ such that  for all $x> 0$,
\begin{eqnarray}\label{ineq5}
  \mathbf{P}\Bigg(  \sum_{i=1}^n\zeta _i \geq x    \Bigg)
   \leq    c_\alpha \exp\Bigg\{ -  \frac{  x^2}{ c_\alpha (u_n +x^{2-\alpha}) }   \Bigg\}.
\end{eqnarray}
\end{lemma}
\noindent
\textit{Proof.} We only give a proof for the case $\alpha \in (0,1).$ For $\alpha=1,$ the proof is similar.
For all $y>0,$ denote $$\eta_{i}(y)=\zeta_{i}\mathbf{1}_{\{\zeta_{i}\leq y\}}.$$
Then  $(\eta_i(y),  \mathcal{F}_{i})_{i=1,...,n}$ is  a sequence of supermartingale differences satisfying $\displaystyle \mathbf{E}(\exp\left\{\lambda \eta_i(y)   \right\} ) < \infty$ for all  $ \lambda \in [0, \infty)$ and all $i$. Define the exponential multiplicative
martingale  $\displaystyle \widetilde{Z}(\lambda )=(\widetilde{Z}_k(\lambda ),\mathcal{F}_k)_{k\geq0},$ where
\[
  \widetilde{Z}_0(\lambda )=1 ,  \quad   \quad \widetilde{Z}_k(\lambda )=\prod_{i=1}^{  k}\frac{\exp\left\{\lambda \eta_i(y)     \right\}}{\mathbf{E}\left(\exp\left\{\lambda \eta_i(y)   \right\} | \mathcal{F}_{i-1} \right)}. \label{C-1}
\]
Then the random variable $\widetilde{Z}_{n}(\lambda ) $  satisfies
$  \int \widetilde{Z}_{n}(\lambda)  d \mathbf{P} = \mathbf{E} \widetilde{Z}_{n}(\lambda) =1.$
Define the  conjugate probability measure
\begin{eqnarray}
d \widetilde{\mathbf{P}}_\lambda =\widetilde{Z}_{ n}(\lambda )d\mathbf{P},   \label{chmeasure3}
\end{eqnarray}
and denote by $\widetilde{\mathbf{E}}_{\lambda}$ the expectation with
respect to $\widetilde{\mathbf{P}}_{\lambda}.$
Notice that $\zeta_i= \eta_i(y) +\zeta_{i}\mathbf{1}_{\{\zeta_{i}> y\}}.$ It is easy to see that   for all $x, y >0$,
\begin{eqnarray}
  \mathbf{P}\Bigg(  \sum_{i=1}^n\zeta _i \geq x  \Bigg)
 \!\!\!&\leq&\!\!\!  \mathbf{P}\left(   \sum_{i=1}^{n}\eta_i(y) \geq x  \right) \ +\  \mathbf{P}\left(   \sum_{i=1}^{n}\zeta_{i}\mathbf{1}_{\{\zeta_{i}> y\}} > 0  \right)  \nonumber\\
&=:&\!\!\!  P_1 + \mathbf{P}\left( \max_{ 1\leq i\leq n}\zeta_{i} > y \right). \label{fmuetoi}
\end{eqnarray}
By the change of measure defined by (\ref{chmeasure3}), we deduce that for all $x, y>0$,
\begin{eqnarray}
  P_1\!\!\! &=&\!\!\! \widetilde{\mathbf{E}}_{\lambda} \Big( \widetilde{Z}_{  n}(\lambda)^{-1}\textbf{1}_{\{\sum_{i=1}^{n}\eta_i(y) \geq x \}} \Big) \nonumber \\
 &=&\!\!\!  \widetilde{\mathbf{E}}_{\lambda}\Big( \exp\Big\{-\lambda \Big(\sum_{i=1}^{n}\eta_i(y) \Big)+ \hat{\Psi}_{n}(\lambda) \Big\} \textbf{1}_{\{\sum_{i=1}^{n}\eta_i(y) \geq x \}} \Big) \nonumber\\
 &\leq&\!\!\!  \widetilde{\mathbf{E}}_{\lambda}\Big( \exp\Big\{-\lambda x + \hat{\Psi}_{n}(\lambda) \Big\} \textbf{1}_{\{\sum_{i=1}^{n}\eta_i(y) \geq x \}} \Big) , \nonumber
\end{eqnarray}
where $\hat{\Psi}_{n}(\lambda)= \sum_{i=1}^n \ln \mathbf{E} \left(e^{\lambda \eta_i(y)  } |
\mathcal{F}_{i-1} \right).$ Notice that
$e^x \leq 1+ x+\frac12 x^2 e^{|x|},\  x \in \mathbf{R}.$ Then, we have for all $\lambda>0$,
\begin{eqnarray*}
\mathbf{E}(e^{\lambda \eta_i(\lambda)}|\mathcal{F}_{i-1})
 \!\!\!&\leq&\!\!\! 1+   \mathbf{E}( \lambda \eta_i(\lambda) |\mathcal{F}_{i-1})+  \frac{\lambda^2}{2} \mathbf{E}( \eta_i(y)^2 \exp\{    |\lambda \eta_{i}(y)|   \}|\mathcal{F}_{i-1}) \\
 &\leq&\!\!\! 1+  \frac{\lambda^2}{2} \mathbf{E}( \eta_i(y)^2 \exp\{\lambda  y^{1-\alpha}   |\eta_{i}(y)|^{\alpha}  \}|\mathcal{F}_{i-1}).
\end{eqnarray*}
 Set $\lambda=c y^{\alpha-1}.$ By  the last inequality and the inequality $\ln(1+t) \leq t$ for all $t\geq 0$,  it is easy to see that for all  $  y>0$,
\begin{eqnarray*}
\hat{\Psi}_{n}(\lambda)\!\!\! &\leq&\!\!\! \sum_{i=1}^n  \ln \bigg( 1+  \frac{\lambda^2}{2} \mathbf{E}( \eta_i(y)^2 \exp\{\lambda  y^{1-\alpha}  |\eta_{i}(y) |^{\alpha}  \}|\mathcal{F}_{i-1}) \bigg)  \\
&\leq&\!\!\! \sum_{i=1}^n \frac{\lambda^2}{2} \mathbf{E}( \eta_i(y)^2 \exp\{\lambda  y^{1-\alpha}  | \eta_{i}(y) |^{\alpha}  \}|\mathcal{F}_{i-1})\\
&\leq&\!\!\! \frac12 c^2 y^{2\alpha-2} u_n.
\end{eqnarray*}
Hence, we get for all $x,  y>0$,
\begin{eqnarray*}
 P_1     \leq  \exp \left\{ - c \,y^{\alpha-1} x +   \frac12 c^2 y^{2\alpha-2} u_n   \right\} .
\end{eqnarray*}
From (\ref{fmuetoi}), it follows that for all $x,  y>0$,
\begin{eqnarray}
  \mathbf{P}\Bigg(  \sum_{i=1}^n\zeta _i \geq x  \Bigg)  \ \leq\
 \exp \left\{ -  c \,y^{\alpha-1} x +   \frac12 c^2 y^{2\alpha-2} u_n   \right\} + \mathbf{P}\left( \max_{ 1\leq i\leq n}\zeta_{i} > y \right).  \label{sdfo}
\end{eqnarray}
By exponential Markov's  inequality, we deduce that for all $y>0$,
\begin{eqnarray}
\mathbf{P}\left( \max_{ 1\leq i\leq n}\zeta_{i} > y \right) \!\!\!&\leq&\!\!\! \sum_{i=1}^n \mathbf{P}\left(  c^{1/\alpha} \zeta_{i} > c^{1/\alpha} y \right) \nonumber\\
 &\leq&\!\!\! \frac{1}{c^{2/\alpha} y^{2}} \exp\{ -c\, y^{\alpha} \} \ \sum_{i=1}^n \mathbf{E} ( \zeta_i^2 \exp\{ c ( \zeta_{i} ^+)^{\alpha} \} ) \nonumber\\
 &\leq&\!\!\!\frac{ u_n}{c^{2/\alpha} y^2}\exp \left\{-c\, y^{\alpha}   \right\}. \label{ft90}
\end{eqnarray}
Taking
\begin{eqnarray*}
y = \left\{ \begin{array}{ll}
\Big(\displaystyle \frac {c\,u_n}{x}\Big)^{1/(1-\alpha)} & \textrm{\ \ \ \ \ if $0< x < (c\,u_n)^{1/(2-\alpha)}$}  \\
\vspace{-0.3cm}\\
x  & \textrm{\ \ \ \ \ if $x \geq (c\,u_n)^{1/(2-\alpha)}$, }
\end{array} \right.
\end{eqnarray*}
from  (\ref{sdfo}) and (\ref{ft90}),    we obtain
\begin{eqnarray}\label{ineq5}
   \mathbf{P}\Bigg(  \sum_{i=1}^n\zeta _i \geq x  \Bigg)
  \leq \left\{ \begin{array}{ll}
\exp\bigg\{\displaystyle \!\! -  \frac{x^2}{2u_n }   \bigg\} +    \frac{ x^{\frac{2}{1-\alpha}}}{u_n^{\frac{1+\alpha}{1-\alpha}}  c^{\frac{2}{\alpha(1+\alpha)} }} \exp\bigg\{ \!\!-c\,\Big(\frac{c\,u_n}{x  } \Big)^{\frac{\alpha}{1-\alpha} } \bigg\}    & \textrm{if\   $0< x < (c\,u_n)^{\frac{1}{2-\alpha}  }$}   \\
\vspace{-0.2cm}\\
\exp \bigg\{\displaystyle\!\! - c\, x^\alpha \Big( 1- \frac{c \, u_n }{2\, x^{2-\alpha}}\Big)  \bigg\} +\frac{ u_n}{c^{2/\alpha} x^2}  \exp \bigg\{\!\! -c\, x^{\alpha}   \bigg\}   & \textrm{if\   $x \geq (c\,u_n)^{\frac{1}{2-\alpha}}$. }
\end{array} \right.
\end{eqnarray}
From (\ref{ineq5}), we get   the following rough bounds for all $x>0,$
\begin{eqnarray}
    \mathbf{P}\Bigg(  \sum_{i=1}^n\zeta _i \geq x  \Bigg)
  \!\!\!&\leq&\!\!\! \left\{ \begin{array}{ll}
c_\alpha \exp\bigg\{\displaystyle -  \frac{\ x^2}{2 c_{1,\alpha} \,u_n }   \bigg\} \ \ \ & \textrm{if\  $0\leq x < (c\,u_n)^{1/(2-\alpha)}$}  \\
\vspace{-0.3cm}  \\
c_\alpha \, \exp \bigg\{\displaystyle  - \frac{c}2 x^\alpha   \bigg\} \ \ \ & \textrm{if\  $x \geq (c\,u_n)^{1/(2-\alpha)}$ }
\end{array} \right. \nonumber  \\
 &\leq & \!\!\! c_\alpha \exp\Bigg\{ -  \frac{  x^2}{ 2 c_\alpha \, (u_n +x^{2-\alpha}) }   \Bigg\}. \nonumber
\end{eqnarray}
This completes the proof of  Lemma \ref{lm21}.
\hfill\qed

In the proof of Theorem \ref{the3.1}, we also need the following normalized Cram\'{e}r-type moderate deviations for martingales; see \cite{FGL13,FS22}.
Let $(\xi _i,\mathcal{F}_i)_{i=0,...,n}$ be a finite sequence of martingale differences.  Set $X_k=\sum_{i=1}^k\xi _i,  k=1,...,n.$
Denote by $\left\langle X\right\rangle $   the quadratic characteristic of the
martingale $X=(X_k,\mathcal{F}_k)_{k=0,...,n},$ that is
\begin{equation}\label{quad}
\left\langle X\right\rangle _0=0,\ \ \ \ \ \ \ \ \left\langle X\right\rangle _k=\sum_{i=1}^k\mathbf{E}(\xi _i^2|\mathcal{F}
_{i-1}),\ \ \ \ \quad k=1,...,n.
\end{equation}
In the sequel we shall use the following conditions:

\begin{description}
\item[(A1)]  There exists  a number $\epsilon_n \in (0, \frac12] $   such that
\[
|\mathbf{E}(\xi_{i}^{k}  | \mathcal{F}_{i-1})| \leq \frac12 k!\epsilon_n^{k-2} \mathbf{E}(\xi_{i}^2 | \mathcal{F}_{i-1}),\ \ \ \ \ \textrm{for all}\ k\geq 2\ \ \textrm{and}\ \ 1\leq i\leq n;
\]
\item[(A2)]  There exist  a number  $ \delta_n\in (0, \frac12]$  and a positive constant $C_1$ such that for all $  x >0 ,$
$$ \mathbf{P}( \left| \left\langle X\right\rangle _n-1\right|  \geq x   )  \leq  C_1 \exp\{  - x^2  \delta_n^{-2} \}.$$
\end{description}

\begin{lemma}\label{th0}
Assume that conditions (A1) and (A2) are satisfied.  Then
the following   inequality holds for all $0 \leq  x   =o( \min\{\epsilon_n^{-1}, $ $  \delta_n ^{-1 } \}),$
\begin{equation}\label{t0ie1}
 \bigg| \ln \frac{\mathbf{P}(X_n/\sqrt{ \langle X\rangle_n } >x)}{1-\Phi \left( x\right)} \bigg|\leq C  \bigg( x^3  (\epsilon_n  + \delta_n)+   (1+  x )  \big(  \delta_n|\ln \delta_n| + \epsilon_n|\ln \epsilon_n| \big) \bigg).
\end{equation}
\end{lemma}

In the proof of Theorem \ref{the3.3}, we   make use of the following lemma.
\begin{lemma}\label{lemma6}
It holds for all $y>0,$
  $$\mathbf{P}\Bigg(  \  \bigg| \sum_{i=0}^{[\eta^{-2}]-1  }  \Big( \eta^{-1 }\big\| ( \theta_{i+1} - \theta_i - \eta g(\theta_i) )^T \nabla \varphi (\theta_i) \big\|_2^2  -  \|\sigma^T \nabla \varphi(\theta_i)\|_2^2      \Big)\bigg|  > y \Bigg) \leq c_1  \exp\bigg\{ -  \frac{y^2}{c_1 (\eta^{-2} +
c\,y ) }   \bigg\}  ,\ \ \ \,$$
  where $c_1$ and  $c $ depend  on $g$ and $\sigma.$
\end{lemma}
\noindent
\textit{Proof.}
Denote by $\psi_{k+1}= \eta^{-1 }\big\| ( \theta_{k+1} - \theta_k - \eta g(\theta_k) )^T \nabla \varphi (\theta_k)\big\|_2^2  -  \|\sigma^T \nabla \varphi(\theta_k)\|_2^2 $. By (\ref{fsdf3}),  the random variable $\psi_{k+1}$ can be rewritten as
$$\psi_{k+1}= \big\| ( \sigma \xi_{k+1})^T \nabla \varphi (\theta_k)\big\|_2^2  -  \|\sigma^T \nabla \varphi(\theta_k)\|_2^2 .$$
 Set $\mathcal{F}_n=\sigma(\theta_0, \xi_{k}, 1 \leq k \leq n  ).$
It is easy to see that $ \mathbf{E} ( \psi_{k+1} |\mathcal{F}_{k})=0  $, and thus $(\psi_i, \mathcal{F}_i)_{i\geq1}$ is a sequence of martingale differences. By Lemma \ref{lm22}, we deduce  that
\begin{eqnarray}
 |\psi_{k+1}| \!\!&\leq&\!\! || \nabla^2 \varphi(\theta_k)  ||  \cdot ||(\sigma \xi_{k+1})(\sigma \xi_{k+1})^{T}-\sigma\sigma^T   || \nonumber \\
 &\leq&\!\! c\, (1+ \| \xi_{k+1}\|_2^2). \nonumber
\end{eqnarray}
The last line implies that there exists a small positive constant $c$ such that $$\mathbf{E}(|  \psi_{k+1}  |^2 \exp\{ c   |  \psi_{k+1}  |    \}   \big|  \mathcal{F}_{k} ) < \infty .$$
Therefore, by Lemma \ref{lm21}, we have for all $  y  \geq c\,\eta^{1/2}$,
 \begin{eqnarray}
&& \mathbf{P}\Bigg(  \  \bigg| \sum_{i=0}^{[\eta^{-2}]-1  }  \Big( \eta^{-1 }\big\| ( \theta_{i+1} - \theta_i - \eta g(\theta_i) )^T \nabla \varphi (\theta_i) \big\|_2^2  -  \|\sigma^T \nabla \varphi(\theta_i)\|_2^2      \Big)\bigg|  > y \Bigg) \nonumber \\
  & &\ \ \ \quad \ \ \ \ \ \quad \ \ \ \   \ \ \ \quad    \ \ \ \ \quad \ \ \ \quad \ \ \ \quad   \ \ \quad \ \ \ \quad \ \ \ \quad   \leq \ \mathbf{P}\bigg(  \Big| \sum_{i=0}^{[\eta^{-2}]-1 }  \psi_{i+1} \Big| >  y\,   \bigg) \nonumber
    \\
  &&\ \ \ \quad \ \ \ \ \quad \ \ \ \ \ \ \quad \ \ \ \     \ \ \ \quad \ \ \ \ \quad \ \ \ \quad \ \ \ \quad   \quad \ \ \ \quad   \leq \  c\, \exp\bigg\{ -  \frac{y^2}{c_1 (\eta^{-2} +
c\,y ) }   \bigg\} . \nonumber
\end{eqnarray}
This completes the proof of Lemma \ref{lemma6}. \hfill\qed

\section{Proof of Theorem \ref{the3.1} }
\setcounter{equation}{0}

Now we are in position to prove Theorem \ref{the3.1}.
Without loss of generality, we assume from now on that $\eta^{-2}$ is an integer. 
From equality (3.1) of Lu, Tan and Xu \cite{LTX22}, we have
\begin{eqnarray*}
\eta^{-1/2}(\Pi_\eta(h)-\pi(h)) \!\!&=&\!\! \mathcal{H}_\eta + \mathcal{R}_\eta,
\end{eqnarray*}
where
$$ \mathcal{H}_\eta=-\eta \sum_{k=0}^{m-1} \langle \nabla \varphi(\theta_k), \sigma \xi_{k+1}  \rangle \ \ \ \ \  \textrm{and} \ \ \  \ \ \  \mathcal{R}_\eta=-\sum_{i=1}^{6}\mathcal{R}_{\eta, i}, $$
with $m=\eta^{-2}$,
\begin{eqnarray*}
\mathcal{R}_{\eta, 1}  \!\! &=&\!\!  \sqrt{\eta}( \varphi(\theta_0)- \varphi(\theta_m) ),\\
\mathcal{R}_{\eta, 2}  \!\! &=&\!\! \frac{\eta^{3/2}}{2} \sum_{k=0}^{m-1}    \langle \nabla^2 \varphi(\theta_k), (\sigma \xi_{k+1})(\sigma \xi_{k+1})^{T}-\sigma\sigma^T  \rangle_{\textrm{HS}}  , \\
\mathcal{R}_{\eta, 3}  \!\! &=&\!\!   \frac{\eta^2}{2} \sum_{k=0}^{m-1}    \langle \nabla^2 \varphi(\theta_k),  g(\theta_k)(\sigma \xi_{k+1})^T\rangle_{\textrm{HS}}  + \langle \nabla^2 \varphi(\theta_k), \sigma\xi_{k+1} g(\theta_k) ^T  \rangle_{\textrm{HS}}  ,\\
\mathcal{R}_{\eta, 4}  \!\! &=&\!\!  \frac{\eta^2}{6} \sum_{k=0}^{m-1}   \int_0^1 \sum_{i_1,i_2,i_3=1}^d \nabla_{i_1,i_2,i_3}^3 \varphi(\theta_k +t \triangle \theta_k) (\sigma \xi_{k+1})_{i_1} (\sigma \xi_{k+1})_{i_2}(\sigma \xi_{k+1})_{i_3}dt,\\
\mathcal{R}_{\eta, 5}  \!\! &=&\!\! \frac{\eta^{5/2}}{2} \sum_{k=0}^{m-1}    \langle \nabla^2 \varphi(\theta_k),  g(\theta_k)g(\theta_k)^T\rangle_{\textrm{HS}} \\
\ && + \  \frac{\eta^{7/2}}{6}  \sum_{k=0}^{m-1}   \int_0^1 \sum_{i_1,i_2,i_3=1}^d \nabla_{i_1,i_2,i_3}^3 \varphi(\theta_k +t \triangle \theta_k) (g(\theta_k))_{i_1} (g(\theta_k))_{i_2}(g(\theta_k))_{i_3}dt , \\
\mathcal{R}_{\eta, 6}  \!\! &=&\!\!  \frac{\eta^{5/2}}{2} \sum_{k=0}^{m-1}   \int_0^1 \sum_{i_1,i_2,i_3=1}^d \bigg[ \nabla_{i_1,i_2,i_3}^3 \varphi(\theta_k +t \triangle \theta_k) (g(\theta_k))_{i_1} (\sigma \xi_{k+1})_{i_2}(\sigma \xi_{k+1})_{i_3} \\
&& \ \ \ \ \ \ \ \ \ \ \ \ \ \ \ \ \ \ \ \ \ \ \ \ \    +  \sqrt{\eta}\, \nabla_{i_1,i_2,i_3}^3\varphi(\theta_k +t \triangle \theta_k) (g(\theta_k))_{i_1} (g(\theta_k))_{i_2}(\sigma \xi_{k+1})_{i_3}   \bigg] dt.
\end{eqnarray*}
Notice that for all $0 \leq x =o(\eta^{-1})$ and  $y = C_0(  \eta^{1/2} x  + (\eta |\ln \eta|)^{1/2})$ with $C_0$ large enough,    we have
\begin{eqnarray*}
\mathbf{P}( W_\eta \geq x) = \mathbf{P}\bigg(\frac{\mathcal{H}_\eta+\mathcal{R}_\eta }{\sqrt{\mathcal{Y}_\eta} }   \geq x  \bigg) \leq \mathbf{P}\bigg( \frac{ \mathcal{H}_\eta}{\sqrt{\mathcal{Y}_\eta} }   \geq x-y  \bigg) + \mathbf{P}\bigg( \frac{\mathcal{R}_\eta }{\sqrt{\mathcal{Y}_\eta} }   \geq y  \bigg)   .
\end{eqnarray*}
Next, we give an estimation for the first term in the r.h.s.\ of the last inequality.
Set $\mathcal{F}_n=\sigma(\theta_0, \xi_{k}, 1 \leq k \leq n  ).$
Then $(-\eta \langle \nabla \varphi(\theta_k), \sigma \xi_{k+1}  \rangle , \mathcal{F}_{k+1})_{k\geq 0}$ is a sequence of martingale differences.
 Since the normal random variable satisfies  the Bernstein condition,  by the boundedness of $|| \nabla \varphi || $ (cf. Lemma \ref{lm22}),   it holds
for all $k\geq 2,$
\[
\Big|\mathbf{E}\big((-\eta \langle \nabla \varphi(\theta_k), \sigma \xi_{k+1}  \rangle )^{k}  \big| \mathcal{F}_{k} \big)\Big| \leq \frac12 k! ( c \|\sigma \|\eta)^{k-2} \mathbf{E}\big( (-\eta   \langle \nabla \varphi(\theta_k), \sigma \xi_{k+1}  \rangle)^2  \big| \mathcal{F}_k \big) \ \
\]
and
$$\langle\mathcal{H}_\eta \rangle_m = \sum_{k=0}^{m-1} \mathbf{E}\big( (-\eta   \langle \nabla \varphi(\theta_k), \sigma \xi_{k+1}  \rangle)^2  \big| \mathcal{F}_k \big)=\mathcal{Y}_\eta.$$
By  Lemma \ref{lm24} with $k=\eta^{-2}$, we have for all $y>0,$
$$ \mathbf{P}\Big( \left| \mathcal{Y}_\eta/\mathbf{E}\mathcal{Y}_\eta-1\right|  \geq y   \Big)  \leq 2 \exp\Big\{  -c\, y^2 \eta^{-2} \Big \} .$$
By Theorem \ref{th0}, we get for  all $0\leq x =o(\eta^{-1})$,
\begin{eqnarray}
\frac{ \mathbf{P}\Big(   \mathcal{H}_\eta / \sqrt{\mathcal{Y}_\eta}     \geq x-y  \Big) }{1- \Phi(x)}\!\!\! &=&\!\!\!\frac{ \mathbf{P}\Big(   \mathcal{H}_\eta / \sqrt{\mathcal{Y}_\eta}     \geq x-y  \Big) }{1- \Phi(x-y)} \ \  \frac{1- \Phi(x-y )}{1- \Phi(x )}  \nonumber \\
\!\! \!&\leq&\!\!\! \exp\bigg\{ c_{1} \bigg(    x^3  \eta  +    (1+x )    \eta |\ln \eta|     \bigg) \bigg\}    \exp\bigg\{ c_2 \,  C_0 x \Big(x\eta^{1/2} + (\eta |\ln \eta|)^{1/2}  \Big)   \bigg\} \nonumber \\
 \!\!\! &\leq&\!\!\!\!  \exp\bigg\{ c_{3} \bigg(    x^3  \eta  +x^2 \eta^{1/2}   +  x(\eta |\ln \eta|)^{1/2}  +   \eta |\ln \eta| \bigg) \bigg\} \nonumber \\
 \!\!\! &\leq&\!\!\!\!  \exp\bigg\{ c_{4} \bigg(    x^3  \eta  +x^2 \eta^{1/2}   + (1+ x) (\eta |\ln \eta|)^{1/2}    \bigg) \bigg\}   .  \label{gfgds}
\end{eqnarray}
Next, we give some estimations for the tail probability $\mathbf{P}(|\mathcal{R}_\eta| \geq y )$ for $0< y =o(\eta^{-1}).$
Clearly, it holds $$\mathbf{P}(|\mathcal{R}_\eta| \geq y) \leq \sum_{i=1}^6\mathbf{P}\big(|\mathcal{R}_{\eta,i}| \geq y/6 \big)=: \sum_{i=1}^6I_i.$$
We now give estimates for $I_i, i=1,2,...,6.$

\emph{\textbf{a)} Control of }$I_1.$
First, by the boundedness of $|\varphi|$ (cf. Lemma \ref{lm22}), we have $|\mathcal{R}_{\eta, 1}| \leq c_1 \, \eta^{1/2}  ,$ and thus for all $y \geq   \eta^{1/2},$
\begin{eqnarray}
I_1 \leq   e^{c_1^2}\, e^{ -   y^2  \eta^{-1}}.
 \end{eqnarray}

\emph{\textbf{b)} Control of }$I_2.$ Denote  $\zeta_{k+1}= \langle \nabla^2 \varphi(\theta_k), (\sigma \xi_{k+1})(\sigma \xi_{k+1})^{T}-\sigma\sigma^T  \rangle_{\textrm{HS}}$.
Then it is easy to see that  $ \mathbf{E} ( \zeta_{k+1} |\mathcal{F}_{k})=0  $ and, by Lemma \ref{lm22}, that
\begin{eqnarray}
 |\zeta_{k+1}| \!\!&\leq&\!\! || \nabla^2 \varphi(\theta_k)   ||  \cdot ||(\sigma \xi_{k+1})(\sigma \xi_{k+1})^{T}-\sigma\sigma^T   || \nonumber \\
 &\leq&\!\! c\, (1+ \| \xi_{k+1}\|_2^2). \nonumber
\end{eqnarray}
The last line implies that there exists a small positive constant $c$ such that $$\mathbf{E}(|   \zeta_{k+1}  |^2 \exp\{ c   |  \zeta_{k+1}  |    \}   \big|  \mathcal{F}_{k} ) < \infty .$$
Therefore, we have for all $  y  \geq c\,\eta^{1/2}$,
 \begin{eqnarray}
I_2  \!\!&\leq&\!\!\!   \mathbf{P}\bigg(  \Big| \sum_{k=0}^{m-1}  \zeta_{k+1} \Big| \geq c_1 \, y\, \eta^{-3/2}  \bigg) \leq   c\, \exp\bigg\{ -  \frac{(c\,y\,\eta^{-3/2})^2}{c_1 (\eta^{-2} +
c\,y\,\eta^{-3/2}) }   \bigg\} \nonumber\\
&\leq&\!\!\!  c_1\, \exp\bigg\{  - \frac{y^2   \eta^{-1}}{c_2\,(1+\eta^{1/2} y )}  \bigg\}.
\end{eqnarray}

\emph{\textbf{c)} Control of }$I_3.$
 The following inequality holds for all $  y  \geq c \,\eta^{1/2}$,
\begin{eqnarray}
I_3   \leq   c_1\, \exp\Big\{  - c_2\, y \eta^{-3/2}  \Big\};
\end{eqnarray}
see   inequality (1.10) in \cite{LTX22a}.

\emph{\textbf{d)} Control of }$I_4.$ It is easy to see that
for all $y>0,$
\begin{eqnarray}
&&I_4:= \mathbf{P}\Big( |\mathcal{R}_{\eta, 4}|  \geq y/6  \Big) \nonumber  \\
  &&\leq   \mathbf{P}\bigg(   \sum_{k=0}^{m-1}   \int_0^1 \sum_{i_1,i_2,i_3=1}^d \nabla_{i_1,i_2,i_3}^3 \varphi(\theta_k +t \triangle \theta_k) (\sigma \xi_{k+1})_{i_1} (\sigma \xi_{k+1})_{i_2}(\sigma \xi_{k+1})_{i_3}dt  \geq    y \eta^{-2}  \bigg) \nonumber  \\
 &&\leq   \mathbf{P}\bigg(   \sum_{k=0}^{m-1} \!   \int_0^1\!\! \sum_{i_1,i_2,i_3=1}^d\!\! \big( \nabla_{i_1,i_2,i_3}^3 \varphi(\theta_k +t \triangle \theta_k) -\nabla_{i_1,i_2,i_3}^3 \varphi(\theta_k )   \big) (\sigma \xi_{k+1})_{i_1} (\sigma \xi_{k+1})_{i_2}(\sigma \xi_{k+1})_{i_3}dt  \geq  \frac12  y \eta^{-2}  \bigg) \nonumber  \\
 &&\ \ \ \ + \  \mathbf{P}\bigg(   \sum_{k=0}^{m-1}   \int_0^1 \!\!\sum_{i_1,i_2,i_3=1}^d \nabla_{i_1,i_2,i_3}^3 \varphi(\theta_k  ) (\sigma \xi_{k+1})_{i_1} (\sigma \xi_{k+1})_{i_2}(\sigma \xi_{k+1})_{i_3}dt  \geq \frac12  y \eta^{-2}  \bigg) \nonumber  \\
 &&\leq   \mathbf{P}\bigg(   \sum_{k=0}^{m-1}\!\! \!  \int_0^1\!\!   \int_0^1\!\!\!  \sum_{i_1,i_2,i_3=1}^d  \nabla_{i_1,i_2,i_3,i_4}^4 \varphi(\theta_k +tt' \triangle \theta_k) (t  \triangle \theta_k)_{i_4}    (\sigma \xi_{k+1})_{i_1} (\sigma \xi_{k+1})_{i_2}(\sigma \xi_{k+1})_{i_3}dt'dt \geq  \frac12  y \eta^{-2}  \bigg) \nonumber  \\
 &&\ \ \ \ + \  \mathbf{P}\bigg(   \sum_{k=0}^{m-1}  \!\! \int_0^1 \sum_{i_1,i_2,i_3=1}^d \nabla_{i_1,i_2,i_3}^3 \varphi(\theta_k  ) (\sigma \xi_{k+1})_{i_1} (\sigma \xi_{k+1})_{i_2}(\sigma \xi_{k+1})_{i_3}dt  \geq \frac12 y \eta^{-2}  \bigg) \nonumber  \\
 & &=: I_{\eta, 4, 1} + I_{\eta, 4, 2}     .  \label{I4}
\end{eqnarray}
We first estimate $I_{\eta, 4, 1}.$
Denote   $\triangle \theta_k=  \theta_{k+1} - \theta_k.$
By the boundedness of   $||\nabla^4\varphi||$ and the fact $\triangle \theta_k=\eta g(\theta_k) + \sqrt{\eta}  \sigma \xi_{k+1} $, we deduce that
for all $y>0,$
\begin{eqnarray*}
I_{\eta, 4, 1}
\!\!\!&\leq&  \!\!\!   \mathbf{P}\bigg(   \sum_{k=0}^{m-1}   \sum_{i_1,i_2,i_3=1}^d   | (  \triangle \theta_k)_{i_4}    (\sigma \xi_{k+1})_{i_1} (\sigma \xi_{k+1})_{i_2}(\sigma \xi_{k+1})_{i_3}|  \geq c_1  y \eta^{-2}  \bigg) \nonumber  \\
 &\leq& \!\!\!  \mathbf{P}\bigg(   \sum_{k=0}^{m-1}   \|    \eta g(\theta_k) + \sqrt{\eta}  \sigma \xi_{k+1} \|_2    \| \sigma \xi_{k+1} \|_2^3 \geq 2c_2  y \eta^{-2}  \bigg) \nonumber  \\
 &\leq&\!\!\! \mathbf{P}\bigg(   \sum_{k=0}^{m-1}   \sqrt{\eta}      \| \sigma \xi_{k+1} \|_2^4 \geq c_2   y \eta^{-2}  \bigg) +   \mathbf{P}\bigg(   \sum_{k=0}^{m-1}   \|    \eta g(\theta_k)  \|_2    \| \sigma \xi_{k+1} \|_2^3 \geq  c_2   y \eta^{-2}  \bigg)     \nonumber \\
 &=:& \!\!\! I_{\eta, 4, 1} ' +  I_{\eta, 4, 1}''.   \nonumber
\end{eqnarray*}
Next, we estimate $I_{\eta, 4, 1} '.$ It is easy to see that
\begin{eqnarray*}
I_{\eta, 4, 1} '  \!\!\!&\leq&\!\!\! \mathbf{P}\bigg(   \sum_{k=0}^{m-1}      \|   \xi_{k+1} \|_2^4 \geq c_2   y \eta^{-5/2}  \bigg)     \nonumber  \\
 &\leq&\!\!\!   \mathbf{P}\bigg( \sum_{k=0}^{m-1} \Big( \|   \xi_{k+1} \|_2^4  -\mathbf{E}  \|   \xi_{k+1} \|_2^4  \Big ) \geq  \Big(c_2 y  \eta^{-5/2}-\sum_{k=0}^{m-1} \mathbf{E}  \|   \xi_{k+1} \|_2^4  \Big )  \bigg)   . \nonumber
\end{eqnarray*}
Clearly, there exists a small positive constant $c$ such that
$$\mathbf{E}\Big(( \|   \xi_{k+1} \|_2^4  -\mathbf{E}  \|   \xi_{k+1} \|_2^4  )^2 \exp\Big\{ c    \big|  \|   \xi_{k+1} \|_2^4  -\mathbf{E}  \|   \xi_{k+1} \|_2^4  \big|^{1/2} \Big\}   \Big|  \mathcal{F}_{k} \Big ) < \infty .$$
Using  Lemma \ref{lm21} with $\alpha=1/2$, we have for all $y\geq c'  \eta^{1/2}$ with $c'$ large enough,
\begin{eqnarray*}
I_{\eta, 4, 1} '
\!\!\! &\leq& \!\!\!   c\, \exp\Bigg\{ -  \frac{ (y  \eta^{-5/2})^2}{c_1 ( \eta^{-2} +(y  \eta^{-5/2})^{3/2}) }   \Bigg\}  \nonumber \\
\!\!\!  &\leq& \!\!\!   c\, \exp\Big\{ -  c_2 \, y^{1/2} \eta^{-5/4}   \Big\} . \nonumber
\end{eqnarray*}
In the sequel, we estimate $I_{\eta, 4, 1} ''$. Using H\"{o}lder's inequality, we get
for all $y>0,$
\begin{eqnarray}
I_{\eta, 4, 1} ''
  \!\!\!&\leq&\!\!\!   \mathbf{P}\bigg( \eta  \Big( \sum_{k=0}^{m-1}   \|     g(\theta_k)  \|_2^2 \Big)^{1/2} \Big( \sum_{k=0}^{m-1}  \| \sigma \xi_{k+1} \|_2^6 \Big)^{1/2}  \geq   c_2  y \eta^{-2}  \bigg)   \nonumber \\
  &\leq&\!\!\!   \mathbf{P}\bigg( \eta  \Big( \sum_{k=0}^{m-1}   \|     g(\theta_k)  \|_2^2 \Big)^{1/2} \Big( \sum_{k=0}^{m-1}  \| \sigma \xi_{k+1} \|_2^6 \Big)^{1/2}  \geq   c_2  y \eta^{-2},\   \eta  \sum_{k=0}^{m-1}   \|     g(\theta_k)  \|_2^2     \geq   C  \,  y^{1/2} \eta^{-5/4}  \bigg)   \nonumber \\
  & & \!\!\! + \ \mathbf{P}\bigg( \eta  \Big( \sum_{k=0}^{m-1}   \|     g(\theta_k)  \|_2^2 \Big)^{1/2} \Big( \sum_{k=0}^{m-1}  \| \sigma \xi_{k+1} \|_2^6 \Big)^{1/2}  \geq   c_2  y \eta^{-2},\   \eta  \sum_{k=0}^{m-1}   \|     g(\theta_k)  \|_2^2     <  C  \,  y^{1/2} \eta^{-5/4}  \bigg)   \nonumber \\
   &\leq&\!\!\!   \mathbf{P}\bigg(   \eta  \sum_{k=0}^{m-1}   \|     g(\theta_k)  \|_2^2     \geq   C  \,  y^{1/2} \eta^{-5/4} \bigg) +  \mathbf{P}\bigg( \eta \Big( Cy^{1/2} \eta^{-9/4} \Big)^{1/2} \Big( \sum_{k=0}^{m-1}  \| \sigma \xi_{k+1} \|_2^6 \Big)^{1/2}  \geq   c_2  y \eta^{-2}  \bigg) \nonumber \\
    &\leq&   \mathbf{P}\bigg(   \eta  \sum_{k=0}^{m-1}   \|     g(\theta_k)  \|_2^2     \geq   C  \,  y^{1/2} \eta^{-5/4} \bigg) +  \mathbf{P}\bigg(  \sum_{k=0}^{m-1}  \| \sigma \xi_{k+1} \|_2^6    \geq   c_2  y^{3/2} \eta^{-15/4}  \bigg) \label{fddf}.
\end{eqnarray}
Next, we give an estimation for the second term in the last inequality.   Denote $T_{k+1}= \| \sigma \xi_{k+1} \|_2^6  -\mathbf{E} \| \sigma \xi_{k+1} \|_2^6.$
It is easy to see  that there exists a small positive constant $c$ such that $$\mathbf{E}\big (  T_{k+1}^{2} \exp\{ c  \ |T_{k+1}|^{1/3} \}   \big|  \mathcal{F}_{k} )  < \infty .$$
Using  Lemma \ref{lm21} with $\alpha=1/3$, we have for all $y\geq c' \eta^{1/2}$ with $c'$ large enough,
\begin{eqnarray*}
  \mathbf{P}\bigg(  \sum_{k=0}^{m-1}  \| \sigma \xi_{k+1} \|_2^6    \geq   c_2   y^{3/2} \eta^{-15/4}  \bigg)
 \!\!\!& =&\!\!\!  \mathbf{P}\bigg( \sum_{k=0}^{m-1} T_{k+1} \geq  \Big(c\,  y^{3/2} \eta^{-15/4}-\sum_{k=0}^{m-1} \mathbf{E}  \| \sigma \xi_{k+1} \|_2^6   \Big)  \bigg) \nonumber  \\
 &\leq&\!\!\! c\, \exp\Bigg\{ -  \frac{ ( y^{3/2} \eta^{-15/4})^2}{c_1 ( \eta^{-2} +( y^{3/2} \eta^{-15/4})^{5/3}) }   \Bigg\}  \nonumber \\
  &\leq& \!\!\!   c\, \exp\Big\{ -  c_2 \, y^{1/2} \eta^{-5/4}   \Big\}.
\end{eqnarray*}
From (\ref{fddf}), by the last inequality and  Lemma \ref{lm23}, we deduce that for all $  y  \geq c' \eta^{1/2}$,
\begin{eqnarray*}
I_{\eta, 4, 1} '' \!\!\!&\leq&\!\!\!   \exp\Big\{ -\gamma  C  \,  y^{1/2} \eta^{-5/4} \Big\} \mathbf{E}\Big[ \exp\Big\{\gamma \eta  \sum_{k=0}^{m-1}   \|     g(\theta_k)  \|_2^2 \Big \} \Big]+ c\, \exp\Big\{ -  c_2 \, y^{1/2} \eta^{-5/4}   \Big\}  \nonumber \\
  &\leq&\!\!\! c_1 \exp\Big\{ - c_2  \,  y^{1/2} \eta^{-5/4} \Big\}   +  c\, \exp\Big\{ -  c_2 \, y^{1/2} \eta^{-5/4}   \Big\}  \nonumber \\
  &\leq&\!\!\! c_3 \exp\Big\{ - c_4 \, y^{1/2} \eta^{-5/4}  \Big\},
\end{eqnarray*}
with  $c' $ and $C$ large enough.
Combining the estimations of $ I'_{\eta, 4, 1}$ and $ I''_{\eta, 4, 2}$, we have the following estimation for $I_{\eta, 4, 1}$: for all $y \geq c'\, \eta^{1/2}$ with $c'$ large enough,
\begin{eqnarray*}
I_{\eta, 4, 1}
 \leq  I'_{\eta, 4, 1} + I''_{\eta, 4, 2}\leq   c_1 \exp\Big\{ - c_2 \, y^{1/2} \eta^{-5/4}  \Big\} .   \nonumber
\end{eqnarray*}

Next, we estimate $I_{\eta, 4, 2}.$
Denote
$$\hbar_{k+1} = \int_0^1 \sum_{i_1,i_2,i_3=1}^d \nabla_{i_1,i_2,i_3}^3 \varphi(\theta_k  ) (\sigma \xi_{k+1})_{i_1} (\sigma \xi_{k+1})_{i_2}(\sigma \xi_{k+1})_{i_3}dt  .$$
Then $(\hbar_{k+1}  , \mathcal{F}_{k+1} )_{k\geq 0} $ is a sequence of martingale differences.
Moreover, by the boundedness of   $||\nabla^3\varphi||$, we have
\begin{eqnarray*}
\big| \hbar_{k+1} \big| \leq c  \sum_{k=0}^{m-1} \|\xi_{k+1}  \|_2^3.
\end{eqnarray*}
Therefore, there exists a small positive constant $c$ such that $\mathbf{E}( |\hbar_{k+1} |^2 \exp\{ c  |\hbar_{k+1} | ^{2/3} \}   \big|  \mathcal{F}_{k} )  < \infty .$
Using Lemma \ref{lm21} with $\alpha =2/3$, we have for all  $y \geq c'\,\eta^{1/2} $ with $c'$ large enough,
\begin{eqnarray*}
I_{\eta, 4, 2}  \!\!\! &\leq&\!\!\!   \mathbf{P}\bigg( \sum_{k=0}^{m-1} \hbar_{k+1}   \geq c y  \eta^{-2}  \bigg)
  =    \mathbf{P}\bigg( \sum_{k=0}^{m-1} \Big(\hbar_{k+1}  -\mathbf{E} \hbar_{k+1} \Big ) \geq \Big(c y  \eta^{-2}-\sum_{k=0}^{m-1} \mathbf{E} \hbar_{k+1}  \Big)  \bigg) \nonumber  \\
  &\leq&\!\!\!    c\, \exp\bigg\{ -  \frac{ (y  \eta^{-2})^2}{c_1 ( \eta^{-2} +(y  \eta^{-2})^{4/3}) }   \bigg\}  \nonumber \\
  &\leq&\!\!\!    c\, \exp\Big\{ -  c_2 \, y^{2/3} \eta^{-4/3}   \Big\}  \nonumber  \\
  &\leq& \!\!\!   c\, \exp\Big\{ -  c_3 \, y^{1/2} \eta^{-5/4}   \Big\} . \nonumber
\end{eqnarray*}
Hence, from (\ref{I4}), we get for all $y \geq c' \eta^{1/2} $ with $c'$ large enough,
\begin{eqnarray*}
I_4  \!\!\!&\leq&\!\!\!  I_{\eta, 4, 1} + I_{\eta, 4, 2}    \leq  c_1 \exp\Big\{  -c_2  \, y^{1/2} \eta^{-5/4}  \Big  \}.
\end{eqnarray*}

\emph{\textbf{e)} Control of }$I_5.$
By the boundedness of $||\nabla^2\varphi||$ and $||\nabla^3\varphi||$ and H\"{o}lder's inequality, we deduce that for all $y>0,$
\begin{eqnarray*}
 I_5 \!\!\! &\leq&\!\!\!  \mathbf{P}\Bigg( \eta^{5/2}  \sum_{k=0}^{m-1}  \| g(\theta_k)\|_2^2  \geq c y  \Bigg) +  \mathbf{P}\Bigg( \eta^{7/2}  \sum_{k=0}^{m-1}  \| g(\theta_k)\|_2^3  \geq c y  \Bigg) \\
 &\leq&\!\!\! \mathbf{P}\Bigg( \eta \sum_{k=0}^{m-1}  \| g(\theta_k)\|_2^2  \geq c y  \eta^{-3/2} \Bigg) + \mathbf{P}\Bigg(\eta^{7/2} m^{1/4}  \Big(\sum_{k=0}^{m-1}  \| g(\theta_k)\|_2^2 \Big )^{3/2}  \geq c y   \Bigg)\\
 &\leq&\!\!\! \mathbf{P}\Bigg( \eta \sum_{k=0}^{m-1}  \| g(\theta_k)\|_2^2  \geq c y  \eta^{-3/2} \Bigg) + \mathbf{P}\Bigg(  \Big(\sum_{k=0}^{m-1}  \| g(\theta_k)\|_2^2 \Big )^{3/2}  \geq c y \eta^{-3}  \Bigg)\\
 &\leq&\!\!\! \mathbf{P}\Bigg( \eta \sum_{k=0}^{m-1}  \| g(\theta_k)\|_2^2  \geq c y  \eta^{-3/2} \Bigg) + \mathbf{P}\Bigg(   \eta\sum_{k=0}^{m-1}  \| g(\theta_k)\|_2^2     \geq c y^{2/3} \eta^{-2}  \Bigg).
\end{eqnarray*}
By   Lemma  \ref{lm23}, we get for all $  y \geq c' \eta^{1/2}  $ with $c'$ large enough,
\begin{eqnarray*}
I_5 \!\!\!&\leq&\!\!\! c \exp\bigg\{ c_1 \eta^{-1}- c_2y  \eta^{-3/2}  \bigg\} +  c \exp\bigg\{ c_1 \eta^{-1}- c_3 y^{2/3} \eta^{-2}   \bigg\}  \\
 & \leq&\!\!\! c_1 \exp\Big\{  -c_2  \, y^{1/2} \eta^{-5/4}  \Big  \}.
\end{eqnarray*}

\emph{\textbf{f)} Control of }$I_6.$
By the boundedness of  $||\nabla^3\varphi||$, one has for all $y> 0,$
\begin{eqnarray}
I_6 \!\!\!&\leq&\!\!\!  \mathbf{P}\Bigg( \eta^{5/2}  \sum_{k=0}^{m-1} \Big( \| g(\theta_k)\|_2 \|  \sigma \xi_{k+1}\|_2 + \sqrt{\eta}   \| g(\theta_k)\|_2^2 \| \sigma \xi_{k+1} \|_2 \Big) \geq c y  \Bigg) \nonumber   \\
 &\leq&\!\!\! \mathbf{P}\Bigg( \sum_{k=0}^{m-1} \| g(\theta_k)\|_2 \|  \sigma \xi_{k+1}\|_2 \geq c y  \eta^{-5/2} \Bigg) + \mathbf{P}\Bigg(\sum_{k=0}^{m-1}   \| g(\theta_k)\|_2^2 \| \sigma \xi_{k+1} \|_2  \geq c y  \eta^{-3}  \Bigg) \nonumber \\
 &=:&\!\!\!  I_{6,1}  + I_{6,2}.
\end{eqnarray}
By H\"{o}lder's inequality, we have for    $C$ large enough and all $y> 0,$
\begin{eqnarray*}
I_{6,1}
\!\!\!&\leq&\!\!\! \mathbf{P}\bigg( \Big(\sum_{k=0}^{m-1}   \| g(\theta_k)\|_2^2 \Big)^{1/2}\Big(\sum_{k=0}^{m-1} \|  \sigma \xi_{k+1}\|_2^2 \Big)^{1/2}     \geq c y  \eta^{- 5/2}  \bigg)  \\
&\leq&\!\!\!   \mathbf{P}\bigg(   \Big( \sum_{k=0}^{m-1}   \|     g(\theta_k)  \|_2^2 \Big)^{1/2} \Big( \sum_{k=0}^{m-1}  \| \sigma \xi_{k+1} \|_2^2 \Big)^{1/2}  \geq   c_2  y \eta^{-5/2},\ \ \ \eta  \sum_{k=0}^{m-1}   \|     g(\theta_k)  \|_2^2     \geq   C  \, y^{1/2}  \eta^{-5/4}  \bigg)   \nonumber \\
  & &\!\!\!  + \ \mathbf{P}\bigg(   \Big( \sum_{k=0}^{m-1}   \|     g(\theta_k)  \|_2^2 \Big)^{1/2} \Big( \sum_{k=0}^{m-1}  \| \sigma \xi_{k+1} \|_2^2 \Big)^{1/2}  \geq   c_2  y \eta^{-5/2},\ \ \ \eta  \sum_{k=0}^{m-1}   \|     g(\theta_k)  \|_2^2     <  C  \,  y^{1/2}  \eta^{-5/4}  \bigg)   \nonumber \\
   &\leq&\!\!\!   \mathbf{P}\bigg(   \eta  \sum_{k=0}^{m-1}   \|     g(\theta_k)  \|_2^2     \geq   C  \, y^{1/2}  \eta^{-5/4}  \bigg) +  \mathbf{P}\bigg( C^{1/2} \Big( \sum_{k=0}^{m-1}  \| \sigma \xi_{k+1} \|_2^2 \Big)^{1/2}  \geq   c_2  y^{3/4} \eta^{-19/8}  \bigg) \nonumber \\
    &\leq&\!\!\!  c_1 \exp\bigg\{-c_2 y^{1/2}  \eta^{-5/4} \bigg\} +  \mathbf{P}\bigg(  \sum_{k=0}^{m-1}  \| \sigma \xi_{k+1} \|_2^2    \geq   c_3  y^{3/2} \eta^{-19/4}  \bigg) .
\end{eqnarray*}
There exists a small positive constant $c$ such that $$\mathbf{E}( \big| \| \sigma \xi_{k+1} \|_2^2  -\mathbf{E}[\| \sigma \xi_{k+1} \|_2^2] \big|^2 \exp\{ c   \big| \| \sigma \xi_{k+1} \|_2^2  -\mathbf{E}[\| \sigma \xi_{k+1} \|_2^2] \big| \}   \big|  \mathcal{F}_{k}  )  < \infty .$$
Using Lemma \ref{lm21} with $\alpha=1$, we have for all $  y \geq   \eta^{1/2} ,$
\begin{eqnarray*}
&& \mathbf{P}\bigg(  \sum_{k=0}^{m-1}  \| \sigma \xi_{k+1} \|_2^2    \geq    c_3  y^{3/2} \eta^{-19/4}   \bigg) \nonumber \\
 &&=   \mathbf{P}\bigg( \sum_{k=0}^{m-1} (\| \sigma \xi_{k+1} \|_2^2  -\mathbf{E}[\| \sigma \xi_{k+1} \|_2^2]  ) \geq (c_1\, y^{3/2} \eta^{-19/4} -\sum_{k=0}^{m-1} \mathbf{E}[\| \sigma \xi_{k+1} \|_2^2]  )  \bigg) \nonumber  \\
  &&\leq   c\, \exp\bigg\{ -  \frac{ (y^{3/2} \eta^{-19/4} )^2}{c_2 ( \eta^{-2} +(y^{3/2} \eta^{-19/4} ) ) }   \bigg\}   \leq  c\, \exp\Big\{ -  c_3 \, y^{3/2 } \eta^{- 19/4}   \Big\}  \nonumber \\
  &&\leq    c\, \exp\Big\{ -  c_3 \, y^{1/2 } \eta^{- 5/4}   \Big\} .
\end{eqnarray*}
Hence,   we get for all $  y \geq   \eta^{1/2} ,$
\begin{eqnarray*}
I_{6,1}   \leq        c_1 \exp\Big\{  -c_2   y^{1/2 } \eta^{- 5/4}  \Big  \}.
\end{eqnarray*}
For $I_{6,2},$  by lemma \ref{lm23},  we have the following estimation for all $  y \geq C \eta^{1/2}  $ with $C$ large enough,
\begin{eqnarray*}
I_{6,2} \!\!\!&=&\!\!\! \mathbf{P}\Bigg(\sum_{k=0}^{m-1}   \| g(\theta_k)\|_2^2 \| \sigma \xi_{k+1} \|_2  \geq c y  \eta^{-3}  \Bigg) \nonumber \\
 &\leq&\!\!\!  \mathbf{P}\bigg(\sum_{k=0}^{m-1}   \| g(\theta_k)\|_2^2  C y^{1/2} \eta^{-3/4}  \geq c y  \eta^{-3}  \bigg)+ \sum_{k=0}^{m-1} \mathbf{P}\bigg( \| \sigma \xi_{k+1} \|_2   \geq C y^{1/2} \eta^{-3/4}   \bigg) \\
 &\leq&\!\!\!  \mathbf{P}\bigg( \eta\sum_{k=0}^{m-1}   \| g(\theta_k)\|_2^2    \geq C   y^{1/2 } \eta^{- 5/4}  \bigg)+ \eta^{-2}\exp\Big\{-C_1 y  \eta^{-3/2}  \Big\} \\
 &\leq&\!\!\!  c_1 \exp\Big\{  - c_2 \, y^{1/2 } \eta^{- 5/4}  \Big\}  + \eta^{-2}\exp\Big\{-C_2 y  \eta^{-3/2}  \Big\}\\
 &\leq&\!\!\!  c_3   \exp\Big\{  - c_4 \, y^{1/2 } \eta^{- 5/4}  \Big\}.
\end{eqnarray*}
Hence, we have for all $  y \geq c' \eta^{1/2}  $ with $c'$ large enough,
\begin{eqnarray}
 I_6 \leq I_{6,1} + I_{6,2}   \leq    c_1  \exp\Big\{  - c_2 \, y^{1/2 } \eta^{- 5/4}  \Big\}.
\end{eqnarray}

Thus, by the estimations of $ I_i, 1\leq i \leq 6,$
we have  for all $  y \geq c' \eta^{1/2}  $ with $c'$ large enough,
\begin{eqnarray}
\mathbf{P}(|\mathcal{R}_\eta| \geq y ) \!\!\!&\leq&\!\!\! c_1 \bigg[ \exp\bigg\{  - \frac{y^2   \eta^{-1}}{c_2\,(1+\eta^{1/2} y )}  \bigg\} + \exp\Big\{  - c_3 \, y^{1/2 } \eta^{- 5/4}  \Big\}   \bigg] \nonumber  . \label{fdsfd}
\end{eqnarray}
Next, we give an estimation for $\mathbf{P}\Big( \frac{\mathcal{R}_\eta }{\sqrt{\mathcal{Y}_\eta} }   \geq y  \Big)$. Clearly, it holds
\begin{eqnarray*}
\mathbf{P}\bigg( \frac{\mathcal{R}_\eta }{\sqrt{\mathcal{Y}_\eta} }   \geq y  \bigg) \!\!\!&\leq&\!\!\! \mathbf{P}\bigg( \frac{\mathcal{R}_\eta }{\sqrt{\mathcal{Y}_\eta} }   \geq y ,  \mathcal{Y}_\eta \leq \mathbf{E}\mathcal{Y}_\eta -\frac12 \mathbf{E}\mathcal{Y}_\eta  \bigg) +   \mathbf{P}\bigg(    \mathcal{Y}_\eta \leq\mathbf{E}\mathcal{Y}_\eta -\frac12 \mathbf{E}\mathcal{Y}_\eta  \bigg)\\
  &\leq& \!\!\! \mathbf{P}\bigg( \frac{\mathcal{R}_\eta }{\sqrt{\mathbf{E} \mathcal{Y}_\eta /2 } }   \geq y   \bigg) +   \mathbf{P}\bigg(    \mathcal{Y}_\eta \leq \mathbf{E}\mathcal{Y}_\eta -\frac12 \mathbf{E}\mathcal{Y}_\eta  \bigg) .
\end{eqnarray*}
By stationarity of $\theta_k$  and $y = C_0 (x\eta^{1/2} + (\eta |\ln \eta|)^{1/2} ),$  we have for all $  y \geq c \, \eta^{1/2}  $ with $c $ large enough,
\begin{eqnarray*}
   \mathbf{P}\bigg( \frac{\mathcal{R}_\eta }{\sqrt{\mathbf{E} \mathcal{Y}_\eta /2 } }   \geq y   \bigg) \!\!\!&\leq&\!\!\! c_1 \bigg[ \exp\bigg\{  - \frac{y^2   \eta^{-1}}{c_2\,(1+\eta^{1/2} y)}  \bigg\}   +  \, \exp\Big\{  - c_3 \, y^{1/2 } \eta^{- 5/4}  \Big\} \bigg] \\
   &\leq&\!\!\! 2 c_1  \bigg[ \exp\bigg\{  - \frac{y^2   \eta^{-1}}{c_2'\,(1+\eta^{1/2} y)}  \bigg\} \mathbf{1}_{\{c \eta^{1/2} \leq y \leq  \eta^{-1/6 }\}} \\
    &&\ \ \ \ \ \  +  \, \exp\Big\{  - c_3' \, y^{1/2 } \eta^{- 5/4}  \Big\}\mathbf{1}_{\{ y > \eta^{-1/6 }\}} \bigg] \\
   &\leq&\!\!\! c_4  \bigg[ \exp\bigg\{  - c_5C_0^2  \big(x   +  \sqrt{|\ln \eta|} \, \big)^2  \bigg\}\mathbf{1}_{\{0 \leq x \leq  \eta^{-2/3 }\}}  \\
   &&\ \ \ \  \ \ \ \ \ + \, \exp\Big\{  - c_6 \, x^{1/2 } \eta^{- 1}  \Big\}\mathbf{1}_{\{x > \eta^{-2/3 }\}} \bigg]
\end{eqnarray*}
and, by Lemma \ref{lm24},
$$\mathbf{P}\Big(    \mathcal{Y}_\eta \leq \mathbf{E}\mathcal{Y}_\eta -\frac12 \mathbf{E}\mathcal{Y}_\eta  \Big) =\mathbf{P}\Big(  \frac12 \mathbf{E}\mathcal{Y}_\eta   \leq \mathbf{E}\mathcal{Y}_\eta - \mathcal{Y}_\eta \Big)  \leq \exp\Big\{-c \, \eta^{-2} \Big\} .$$
Hence, for all $0<  x \leq \eta^{-2/3 } $,
\begin{eqnarray}
\mathbf{P}\bigg( \frac{\mathcal{R}_\eta }{\sqrt{\mathcal{Y}_\eta} }   \geq y  \bigg)   \leq  c    \exp\bigg\{  - c_5C_0^2  \big(x   +  \sqrt{|\ln \eta|} \, \big)^2  \bigg\}     . \label{hgfhd}
\end{eqnarray}
Take $C_0$ such that $c_5C_0^2 \geq 4.$ Combining the inequalities (\ref{gfgds}) and (\ref{hgfhd}) together, we get
 for   all  $0\leq x   \leq \eta^{-2/3 } $,
\begin{eqnarray}
\frac{ \mathbf{P}( W_\eta \geq x) }{1- \Phi(x)} \!\!\!&\leq&\!\!\! \frac{ \mathbf{P}\big(   \mathcal{H}_\eta / \sqrt{\mathcal{Y}_\eta}     \geq x-y  \big) }{1- \Phi(x)}  + \frac{\mathbf{P}\big(  \mathcal{R}_\eta /\sqrt{\mathcal{Y}_\eta}     \geq y  \big)}{1- \Phi(x)}  \nonumber \\
  &\leq&\!\!\! \exp\bigg\{ c_{2} \bigg(    x^3  \eta  +x^2 \eta^{1/2}   +  x (\eta |\ln \eta|)^{1/2}  +     \eta |\ln \eta| \bigg) \bigg\}  \nonumber \\
    && + \frac{c}{1- \Phi(x) }  \exp\bigg\{  -  c_5C_0^2  \big(x   +  \sqrt{|\ln \eta|} \, \big)^2  \bigg\}  \nonumber \\
    &\leq&\!\!\! \exp\bigg\{ c_{3} \bigg(    x^3  \eta  +x^2 \eta^{1/2}   +  (1+x) (\eta |\ln \eta|)^{1/2}   \bigg) \bigg\}    \label{deineq01} .
\end{eqnarray}
Similarly,  we have for all $0 \leq x \leq \eta^{-2/3} $,
\begin{eqnarray}
\frac{ \mathbf{P}( W_\eta \geq x) }{1- \Phi(x)}\!\!\! &\geq&\!\!\! \frac{ \mathbf{P}\big(   \mathcal{H}_\eta / \sqrt{\mathcal{Y}_\eta}     \geq x+y  \big) }{1- \Phi(x)}  - \frac{\mathbf{P}\big(  \mathcal{R}_\eta /\sqrt{\mathcal{Y}_\eta}    \leq - y  \big)}{1- \Phi(x)}  \nonumber \\
  &\geq&\!\!\! \exp\bigg\{ -c_{2}\bigg(    x^3  \eta  +x^2 \eta^{1/2}   + (1+ x) (\eta |\ln \eta|)^{1/2}  \bigg) \bigg\}    \nonumber\\
  &&\  -  \frac{c}{1- \Phi(x) }  \exp\bigg\{  -  c_5C_0^2  \big(x   +  \sqrt{|\ln \eta|} \, \big)^2  \bigg\}  \nonumber\\
    &\geq&\!\!\! \exp\bigg\{- c_{3} \bigg(    x^3  \eta  +x^2 \eta^{1/2}   + (1+ x) (\eta |\ln \eta|)^{1/2}   \bigg)   \bigg\}    .\label{deineq02}
\end{eqnarray}
Combining the inequalities (\ref{deineq01}) and (\ref{deineq02}) together, we obtain  the first desired inequality for all $0 \leq x \leq \eta^{-2/3} $.
For the case $ \eta^{-2/3} < x \leq \eta^{-3/4},$  the assertion of Theorem \ref{the3.1} follows by a similar argument by taking  $y=C_0 x^{2}\eta $, instead of  $y = C_0 (x\eta^{1/2} + (\eta |\ln \eta|)^{1/2}  ) $,  and accordingly in the subsequent statements.

 The result for $-W_\eta$  follows by the first inequality applying to $-W_\eta$.  This completes the proof of
 Theorem \ref{the3.1}.

 \section{Proof of Corollary \ref{coro1}}
\setcounter{equation}{0}
 It is easy to see that
\begin{eqnarray}
  \sup_{x \in \mathbf{R}}\Big|\mathbf{P}( W_\eta   \leq  x  ) - \Phi \left(  x\right)\Big|
\!\!\!&\leq&\!\!\! \sup_{   x > \eta^{-1/8}} \big|\mathbf{P}(W_\eta  \leq  x  ) - \Phi \left(  x\right) \big|    + \sup_{  0 \leq x \leq  \eta^{-1/8} } \big|\mathbf{P}(W_\eta   \leq  x  ) - \Phi \left(  x\right)\big| \nonumber\\
&+ &\!\!\!  \sup_{ - \eta^{-1/8} \leq x \leq 0 } \big|\mathbf{P}(W_\eta   \leq  x  ) - \Phi \left(  x\right) \big|  +\sup_{   x < - \eta^{-1/8}} \big|\mathbf{P}(W_\eta   \leq  x  ) - \Phi \left(  x\right) \big| \nonumber  \\
&=:\!\!\! & H_1 + H_2+H_3+H_4.   \label{indeq0d10}
\end{eqnarray}
It is known that
\begin{equation}
\frac 1{\sqrt{ 2 \pi}(1+ \lambda)}\ e^{- \lambda^2/2} \leq 1-\Phi \left(  \lambda\right)
\leq
\frac 1{\sqrt{  \pi}(1+ \lambda)}\ e^{- \lambda^2/2},\ \ \ \ \lambda\geq 0,
\label{f39}
\end{equation}
see \cite{GH00}.
By Theorem \ref{the3.1}, we deduce that
\begin{eqnarray*}
H_1 & \leq &\!\!\!    \sup_{   x > \eta^{-1/8} }  \mathbf{P}\big( W_\eta  > x   \big) + \sup_{   x > \eta^{-1/8} }  \big(1 -  \Phi \left( x\right) \big)
   \ \leq \  \mathbf{P}\big( W_\eta > \eta^{-1/8}   \big) +    \big(1 -  \Phi  ( \eta^{-1/8}  ) \big) \\
   &\leq& \!\!\!   \big(1 -  \Phi  ( \eta^{-1/8}  ) \big)e^c   +    \exp\Big\{ -\frac12  \eta^{-1/4}  \Big \}  \\
 & \leq &\!\!\!  c_1 (\eta  |\ln \eta|)^{1/2}
\end{eqnarray*}
and
\begin{eqnarray*}
H_4   &\leq &\!\!\! \sup_{   x <- \eta^{-1/8} }  \mathbf{P}\big( W_\eta    \leq x  \big) + \sup_{   x <- \eta^{-1/8} }     \Phi \left( x\right)
   \ \leq  \    \mathbf{P}\big( W_\eta \leq - \, \eta^{-1/8}  \big) +     \Phi  ( -\, \eta^{-1/8}   )  \\
 &\leq &\!\!\!   \Phi (-\, \eta^{-1/8} )  e^c   +    \exp\Big\{ -\frac12 \eta^{-1/4} \Big\} \\
 & \leq &\!\!\!  c_2   (\eta  |\ln \eta|)^{1/2} .
\end{eqnarray*}
Using Theorem \ref{the3.1} and the inequality $|e^x-1|\leq |x|e^{|x|},$ we get
\begin{eqnarray*}
H_2 \!\!\!&=&\!\!\! \sup_{0\leq   x \leq \eta^{-1/8} } \Big|\mathbf{P}\big(W_\eta  > x  \big)- \big(1 -  \Phi \left( x\right) \big) \Big| \nonumber \\
 &\leq&\!\!\!  \sup_{0\leq   x \leq \eta^{-1/8} } c_{ 1} \Big(1-\Phi(x) \Big) \bigg(   x^3  \eta  +x^2 \eta^{1/2}   + (1+ x) (\eta  |\ln \eta|)^{1/2}      \bigg)  \nonumber\\
 &\leq &\!\!\! c_{3}  (\eta  |\ln \eta|)^{1/2}
 \end{eqnarray*}
 and
 \begin{eqnarray*}
H_3 &=&\!\!\! \sup_{ - \eta^{-1/8} \leq x \leq 0 } \big|\mathbf{P}\big( W_\eta  \leq x  \big) -  \Phi \left( x\right) \big|  \nonumber \\
   &\leq&\!\!\!   \sup_{ - \eta^{-1/8} \leq x \leq 0 } c_{ 1}  \Phi(x)  \bigg(   |x|^3  \eta  +x^2 \eta^{1/2}   + (1+ |x|) (\eta  |\ln \eta|)^{1/2}      \bigg)  \nonumber\\
 &\leq &\!\!\! c_{ 4} \, (\eta  |\ln \eta|)^{1/2} .
 \end{eqnarray*}
Applying  the bounds  of $H_1, H_2, H_3$ and $H_4$ to (\ref{indeq0d10}),   we obtain the desired  inequality of Corollary \ref{coro1}.

\section{Proof of Theorem \ref{the3.3} }
\setcounter{equation}{0}
Assume that $ \varepsilon_x \in (0, 1/2].$
It is easy to see that for all $x\geq 0,$
\begin{eqnarray}
 \mathbf{P}\Big( S_\eta > x   \Big) \!\!\!& =& \!\!\! \mathbf{P}\bigg(  \eta^{-1/2}   (  \Pi_\eta(h) - \pi (h) ) > x \sqrt{ \mathcal{Y}_\eta},\  \mathcal{Y}_\eta \geq (1- \varepsilon_x)\mathcal{V}_\eta \bigg) \nonumber  \\
 && + \ \mathbf{P}\bigg(  \eta^{-1/2}   (  \Pi_\eta(h) - \pi (h) )  > x \sqrt{ \mathcal{Y}_\eta},\   \mathcal{Y}_\eta <(1- \varepsilon_x)\mathcal{V}_\eta  \bigg) \nonumber  \\
 &\leq&\!\!\! \mathbf{P}\bigg( W_\eta \geq x \sqrt{ 1- \varepsilon_x  \ }   \bigg)
  + \ \mathbf{P}\bigg(  \mathcal{Y}_\eta - \mathcal{V}_\eta   < - \varepsilon_x \mathcal{V}_\eta , \ \mathcal{V}_\eta   \geq \frac12 \mathbf{E}\mathcal{V}_\eta \bigg) \nonumber  \\
   && + \  \mathbf{P}\bigg(  \mathcal{Y}_\eta - \mathcal{V}_\eta   < - \varepsilon_x \mathcal{V}_\eta , \  \mathcal{V}_\eta   < \frac12 \mathbf{E}\mathcal{V}_\eta \bigg)  \nonumber  \\
  &\leq&\!\!\! \mathbf{P}\bigg( W_\eta \geq x \sqrt{ 1- \varepsilon_x  \ }   \bigg)
  + \ \mathbf{P}\bigg(  \mathcal{Y}_\eta - \mathcal{V}_\eta   < -\frac12  \varepsilon_x \mathbf{E}\mathcal{V}_\eta    \bigg) + \mathbf{P}\bigg(  \mathcal{V}_\eta  - \mathbf{E}\mathcal{V}_\eta   < -\frac12 \mathbf{E}\mathcal{V}_\eta    \bigg) \nonumber  \\
 &=:&\!\!\! P_1 +P_2 +P_3.    \label{a15}
\end{eqnarray}
By Theorem  \ref{the3.1}, we have for all   $0 \leq x \leq \eta^{-3/4} $,
\begin{eqnarray}\label{fsds01}
 P_1  \!\!\! & \leq&\!\!\! \Big(1-\Phi(x\sqrt{1- \varepsilon_x })\Big)\exp\bigg\{    c  \Big(      x^3  \eta  +x^2 \eta^{1/2}   + (1+ x) (\eta  |\ln \eta|)^{1/2}      \Big)  \bigg\} \nonumber  \\
 & \leq&\!\!\! \Big(1-\Phi(x)\Big)\exp\bigg\{    c  \Big(x \varepsilon_x   +    x^3  \eta  +x^2 \eta^{1/2}   + (1+ x) (\eta  |\ln \eta|)^{1/2}   \Big)  \bigg\} .
\end{eqnarray}
Using Lemma  \ref{lm24}, we get for all $x\geq 0,$
\begin{eqnarray}\label{fsds02}
 P_2   \  \leq \ 2\, \exp \Big\{- c \,  \varepsilon_x ^2 \eta^{-2} \Big\}.
\end{eqnarray}
By Lemma \ref{lemma6}, we deduce that for all $x \geq 0,$
\begin{eqnarray}\label{fsds03}
 P_3 \ \leq \ c_1  \exp\Big\{ - c \,  \eta^{-2}  \Big\}  .
\end{eqnarray}
Taking $\varepsilon_x= c_0 x\eta  + \eta  ^{1/2}  $ with $c_0$ large enough, by (\ref{a15})-(\ref{fsds03}),
we deduce that for all  $0 \leq x \leq \eta^{-3/4} $,
\begin{eqnarray}
 \mathbf{P}\Big( S_\eta > x   \Big)
 \!\!\!&\leq&\!\!\!\Big(1-\Phi(x)\Big)  \exp\bigg\{    c  \Big(   x^3  \eta  +x^2 \eta^{1/2}   + (1+ x) (\eta  |\ln \eta|)^{1/2}   \Big)  \bigg\} \nonumber  \\
 && + \  2\, \exp \bigg\{- c \, \Big( c_0^2 \, x^2     +     \eta^{-1}   \Big )    \bigg\} \ + \ c_1  \exp\bigg\{ - c \,  \eta^{-2}  \bigg\} . \nonumber
\end{eqnarray}
Applying (\ref{f39}) to the last inequality, we obtain for all $ 0 \leq x \leq \eta^{-3/4} $,
\begin{eqnarray}\label{dfschj}
\mathbf{P}\Big( S_\eta > x   \Big)
 \ \leq \ \Big(1-\Phi(x)\Big)  \exp\bigg\{    c  \Big(   x^3  \eta  +x^2 \eta^{1/2}   + (1+ x) (\eta  |\ln \eta|)^{1/2}   \Big)  \bigg\} ,
\end{eqnarray}
which gives the  upper bound  for the tail probability $\mathbf{P}\big( S_\eta > x   \big), x\geq 0.$
Notice that for all $x\geq 0,$
\begin{eqnarray}
 \mathbf{P}\Big( S_\eta > x   \Big) & \geq& \ \mathbf{P}\bigg(  \eta^{-1/2}   (  \Pi_\eta(h) - \pi (h) )  > x \sqrt{ \mathcal{Y}_\eta},\   \mathcal{Y}_\eta <(1+ \varepsilon_x)\mathcal{V}_\eta  \bigg) \nonumber  \\
 &\geq&\mathbf{P}\bigg( W_\eta \geq x \sqrt{ 1+ \varepsilon_x  \ }   \bigg)
  - \ \mathbf{P}\bigg(  \mathcal{Y}_\eta - \mathcal{V}_\eta   \geq \varepsilon_x \mathcal{V}_\eta , \ \mathcal{V}_\eta    \leq \frac12 \mathbf{E}\mathcal{V}_\eta \bigg) \nonumber  \\
   && - \  \mathbf{P}\bigg(  \mathcal{Y}_\eta - \mathcal{V}_\eta   \geq\varepsilon_x \mathcal{V}_\eta , \  \mathcal{V}_\eta   > \frac12 \mathbf{E}\mathcal{V}_\eta \bigg)  \nonumber  \\
  &\geq&\mathbf{P}\bigg( W_\eta \geq x \sqrt{ 1+ \varepsilon_x  \ }   \bigg)
  - \ \mathbf{P}\bigg(  \mathcal{V}_\eta -  \mathbf{E} \mathcal{V}_\eta  \leq -  \frac12 \mathbf{E}\mathcal{V}_\eta    \bigg)
   -\mathbf{P}\bigg(  \mathcal{Y}_\eta - \mathcal{V}_\eta   \geq   \frac12\varepsilon_x \mathbf{E}\mathcal{V}_\eta    \bigg) \nonumber  \\
 &=:& P_4 -P_5  - P_6.    \label{as1f5}
\end{eqnarray}
By Theorem  \ref{the3.1}, we have for all   $0 \leq x \leq \eta^{-3/4} $,
\begin{eqnarray}\label{adfsds01}
 P_4  \!\!\! & \geq&\!\!\! \Big(1-\Phi(x\sqrt{1+ \varepsilon_x })\Big)\exp\bigg\{   - c  \Big(      x^3  \eta  +x^2 \eta^{1/2}   + (1+ x) (\eta  |\ln \eta|)^{1/2}      \Big)  \bigg\} \nonumber  \\
 & \geq&\!\!\! \Big(1-\Phi(x)\Big)\exp\bigg\{  -  c  \Big(x \varepsilon_x   +    x^3  \eta  +x^2 \eta^{1/2}   + (1+ x) (\eta  |\ln \eta|)^{1/2}   \Big)  \bigg\} .
\end{eqnarray}
Using Lemma  \ref{lm24} with $k=\eta^{-2}$, we get for all $x\geq 0,$
\begin{eqnarray}\label{adfsds02}
 P_5  \ \leq \ 2\, \exp \Big\{- c \,   \eta^{-2} \Big\}.
\end{eqnarray}
By Lemma \ref{lemma6}, we deduce that for all $x \geq 0,$
\begin{eqnarray}  \label{adfsds03}
 P_6 \ \leq \  c_1  \exp\bigg\{ -  \frac{(\varepsilon_x \eta^{-2} )^2}{c_1 (\eta^{-2} +
c\,\varepsilon_x \eta^{-2}  ) }   \bigg\} \ \leq \ c_1  \exp\bigg\{ -   c_2 \, \varepsilon_x  ^2 \eta^{-2}  \bigg\}.
\end{eqnarray}
Taking $\varepsilon_x= c_0 x\eta  + \eta  ^{1/2} $ with $c_0$ large enough, by (\ref{as1f5})-(\ref{adfsds03}),
we deduce that for all  $0 \leq x \leq \eta^{-3/4} $,
\begin{eqnarray}
 \mathbf{P}\Big( S_\eta > x   \Big)
\!\!\! &\geq&\!\!\! \Big(1-\Phi(x)\Big)  \exp\bigg\{  -  c  \Big(   x^3  \eta  +x^2 \eta^{1/2}   + (1+ x) (\eta  |\ln \eta|)^{1/2}   \Big)  \bigg\} \nonumber  \\
 && - \, 2\, \exp \bigg\{- c \,    \eta^{-2}        \bigg\} - c_1  \exp\bigg\{ - c \, \Big( c_0^2\, x^2     +   \eta^{-1}     \Big ) \bigg\} . \nonumber
\end{eqnarray}
Applying (\ref{f39}) to the last inequality, we obtain for all $ 0 \leq x \leq \eta^{-3/4} $,
\begin{eqnarray}
\mathbf{P}\Big( S_\eta > x   \Big)
 \ \geq \ \Big(1-\Phi(x)\Big)  \exp\bigg\{ -   c  \Big(   x^3  \eta  +x^2 \eta^{1/2}   + (1+ x) (\eta  |\ln \eta|)^{1/2}   \Big)  \bigg\} ,
\end{eqnarray}
which gives the lower bound  for the tail probability $\mathbf{P}\big( S_\eta > x   \big), x\geq 0.$
The proof for $- S_\eta$  follows by a similar argument. This completes the proof of Theorem \ref{the3.3}.

\section*{Acknowledgements}
The work has been partially supported by the National Natural Science Foundation of China (Grant no.\ 11971063).

\end{document}